\documentclass[12pt]{article}

\usepackage{amsfonts,amsmath,amssymb}
\usepackage{mathrsfs,mathtools,url}
\usepackage{latexsym}
\usepackage{tikz,color}
\usepackage[T1]{fontenc}
\usepackage{array}
\usepackage{xspace}
\usepackage{cite}

\newtheorem{theorem}{Theorem}[section]
\newtheorem{definition}[theorem]{Definition}
\newtheorem{lemma}[theorem]{Lemma}
\newtheorem{proposition}[theorem]{Proposition}
\newtheorem{observation}[theorem]{Observation}

\newtheorem{corollary}[theorem]{Corollary}
\newtheorem{remark}[theorem]{Remark}

\newtheorem{claim}{Claim}

\DeclareMathOperator {\diam} {diam}

\let\deg\relax
\DeclareMathOperator {\deg} {deg}
\DeclareMathOperator {\mut} {\mu_{\rm t}}
\DeclareMathOperator {\muo} {\mu_{\rm o}}
\DeclareMathOperator {\mud} {\mu_{\rm d}}
\DeclareMathOperator {\bp} {bp}
\DeclareMathOperator {\FS} {FS}
\DeclareMathOperator {\name} {{\cal V}}
\DeclareMathOperator {\nameo} {{\cal V}_{\rm o}}
\DeclareMathOperator {\named} {{\cal V}_{\rm d}}
\DeclareMathOperator {\namet} {{\cal V}_{\rm t}}

\newcommand{\MV}{mutual-visibility}

\newcommand{\proof}{\noindent{\bf Proof.\ }}
\newcommand{\qed}{\hfill $\square$ \bigskip}
\newcommand{\smallqed}{{\tiny ($\Box$)}}

\newcommand{\cG}{{\cal G}}
\newcommand{\cF}{{\cal F}}
\newcommand{\cP}{{\cal P}}
\newcommand{\sandi}[1]{\textcolor{green!60!black}{#1}}

\title{Visibility polynomials, dual visibility spectrum, and characterization of total mutual-visibility sets}

\author{Csilla Bujt\'{a}s$^{a, b}$, Sandi Klav\v{z}ar$^{a, b, c}$, Jing Tian$\/^{d,b}$ \\\\
$^{a}$ \small Faculty of Mathematics and Physics, University of Ljubljana, Slovenia\\
\small {\tt csilla.bujtas@fmf.uni-lj.si\quad ORCID: 0000-0002-0511-5291}\\
\small {\tt sandi.klavzar@fmf.uni-lj.si\quad ORCID: 0000-0002-1556-4744}\\
$^{b}$ \small Institute of Mathematics, Physics and Mechanics, Ljubljana, Slovenia \\
$^{c}$ \small Faculty of Natural Sciences and Mathematics, University of Maribor, Slovenia\\
$^{d}$ \small School of Science, Zhejiang University of Science and Technology, \\
\small Hangzhou, Zhejiang 310023, PR China\\
\small {\tt jingtian526@126.com \quad ORCID: 0000-0002-1578-4798}\\
}

\date{\today}

\begin{document}

\maketitle

\begin{abstract}
Mutual-visibility sets were motivated by visibility in distributed systems and social networks, and intertwine with several classical mathematical areas. Monotone properties of the variety of mutual-visibility sets, and restrictions of such sets to convex and isometric subgraphs are studied. Dual mutual-visibility sets are shown to be intrinsically different from other types of mutual-visibility sets. It is proved that for every finite subset $Z$ of positive integers there exists a graph $G$ that has a dual mutual-visibility set of size $i$ if and only if $i\in Z\cup \{0\}$, while for the other types of mutual-visibility such a set consists of consecutive integers. Visibility polynomials are introduced and their properties derived. As a surprise, every polynomial with nonnegative integer coefficients and with a constant term $1$ is a dual visibility polynomial of some graph. Characterizations are given for total mutual-visibility sets, for graphs with total mutual-visibility number $1$, and for sets which are not total mutual-visibility sets, yet every proper subset is such. Along the way an earlier result from the literature is corrected.
\end{abstract}


\medskip\noindent
\textbf{AMS Math.\ Subj.\ Class.\ (2020)}: 05C12, 05C31, 05C69


\section{Introduction}

Let $G = (V(G), E(G))$ be a graph and $X\subseteq V(G)$. Then vertices $x$ and $y$ of $G$ are {\em $X$-visible}, if there exists a shortest $x,y$-path $P$ such that no internal vertex of $P$ belongs to $X$. The set $X$ is a \emph{mutual-visibility set} if any two vertices from $X$ are $X$-visible, while $X$ is a \emph{total mutual-visibility set} if any two vertices from $V(G)$ are $X$-visible. Let $\overline{X} = V(G)\setminus X$. Then $X$ is a \emph{dual mutual-visibility set} if any two vertices from $X$ and any two vertices from $\overline{X}$ are $X$-visible. Finally, $X$ is an \emph{outer mutual-visibility set} if any two vertices from $X$ are $X$-visible, and any two vertices $x\in X$, $y\in \overline{X}$ are $X$-visible. The cardinality of a largest  mutual-visibility set (resp.\ total/dual/outer mutual-visibility set) is the \emph{mutual-visibility number} (resp.\ \emph{total/dual/outer mutual-visibility number}) $\mu(G)$ (resp.\ $\mut(G)$, $\mud(G)$, $\muo(G)$) of $G$. A mutual-visibility set of cardinality $\mu(G)$ is called a {\em $\mu$-set}. We have analogous meaning for $\mut$-sets, $\mud$-sets, and $\muo$-sets. The key definitions are summarized in Table~\ref{tab:example}, where for arbitrary vertices $x,y\in V(G)$, we denote by ``$+$'' if $x$ and $y$ are required to be $X$-visible, and by ``$-$'' if it is not required. 

\begin{table}[ht!]
\centering
\begin{tabular}{|>{\centering\arraybackslash}m{5cm}|>{\centering\arraybackslash}m{1.7cm}|>{\centering\arraybackslash}m{1.8cm}|>{\centering\arraybackslash}m{1.8cm}|>
{\centering\arraybackslash}m{2cm}|}
\hline \vspace{25pt}
 \raisebox{0pt}[10pt][10pt] & $x \in X$, $y \in X$
 & $x \in X$, $y \in \overline{X}$ 
 & $x \in \overline{X}$, $y \in \overline{X}$ 
 & maximum cardinality 
 \\ \hline 

 mutual-visibility & $+$ & $-$ & $-$ & $\mu$ \\ \hline
 
 dual mutual-visibility & $+$ & $-$ & $+$ & $\mu_{\rm d}$ \\ \hline
 
outer  mutual-visibility & $+$ & $+$ & $-$ & $\mu_{\rm o}$ \\ \hline
 
 total mutual-visibility & $+$ & $+$ & $+$ & $\mu_{\rm t}$ \\ \hline
\end{tabular}
\caption{Varity of visibility definitions}
\label{tab:example}
\end{table}

Mutual-visibility sets were introduced by Di Stefano in~\cite{distefano-2022} motivated by mutual visibility in distributed computing and social networks. Although the motivation came from theoretical computer science, it is a graph theoretical concept. It needs to be said that the term mutual-visibility is also used in other contexts, for instance in robotics, where the mutual visibility problem asks for a distributed algorithm that repositions robots to a configuration where they all can see each other, cf.~\cite{adhikary-2022}. 

The graph theoretic mutual-visibility has received a lot of interest and was investigated in a series of papers~\cite{axenovich-2024a+, boruzanli-2024, BresarYero-2024, cicerone-2025, cicerone-2023a, cicerone-2023+, cicerone-2023, cicerone-2024b, korze-2024, korze-2024+, roy-2025}. In addition to being an interesting concept, the fact that the topic is intertwined with several other areas has also contributed to the interest. These include the Zarankiewicz problem~\cite{cicerone-2023}, Tur\'an type problems on graphs and hypergraphs~\cite{boruzanli-2024, Bujtas, cicerone-2024b}, and a close relationship with the Bollob\'as-Wessel theorem~\cite{bollobas-1967, wessel-1966} as established in~\cite{BresarYero-2024}. Also, Axenovich and Liu~\cite{axenovich-2024a+} proved that $\mu(Q_n) \ge 0.186\cdot 2^n$ by using a recent  breakthrough result on daisy-free hypergraphs due to Ellis, Ivan, and Leader~\cite{ellis-2024+}.

The investigations from~\cite{cicerone-2023+} raised the need to introduce the total mutual-visibility which was in turn studied in~\cite{axenovich-2024a+, boruzanli-2024, Bujtas, kuziak-2023, tian-2024}. The remaining two types of visibility were introduced in~\cite{CiDiDrHeKlYe-2023} and further considered in~\cite{BresarYero-2024, korze-2024+, roy-2025}. 

In this paper we first consider monotone properties of the variety of mutual-visibility sets, and restrictions of such sets to convex and isometric subgraphs. Along the way an earlier result from the literature is corrected. In Section~\ref{sec:polynomials} we introduce visibility polynomials, show some examples, and derive some properties of these polynomials. Since it is observed in Section~\ref{sec:monotonicity} that dual mutual-visibility sets are intrinsically different from other types of mutual-visibility sets, we introduce in Section~\ref{sec:dual-spectrum} the dual visibility spectrum as the counting vector of dual mutual-visibility sets of different sizes. The main result of the section shows that the nonnegative entries can be arbitrarily prescribed and a graph with this visibility spectrum exists. In other words, every polynomial with nonnegative integer coefficients and with a constant term $1$ is a dual visibility polynomial of some graph. In the final section we consider total mutual-visibility sets. We give a general characterization, describe graphs $G$ with $\mut(G) = 1$, and characterize sets which are not total mutual-visibility sets, yet every proper subset is such. 

\medskip
In the rest of the introduction we give additional definitions needed. If $G$ is a graph and $v\in V(G)$, then $N_G(v)$ denotes the set of vertices adjacent to $v$. The {\em degree} $\deg_G(v)$ of $v$ is $|N_G(v)|$. 

For vertices $u$ and $v$ of $G$, the length of a shortest $u,v$-path is the {\em distance} between $u$ and $v$ and denoted by $d_G(u,v)$. A subgraph $H$ of $G$ is {\em isometric}, if $d_H(u,v) = d_G(u,v)$ for every two vertices $u$ and $v$ of $H$. Further, $H$ is {\em convex}, if for every two vertices of $H$, all shortest $u,v$-paths belong to $H$. A graph $G$ is {\em geodetic} if the shortest path between each pair of vertices is unique, cf.~\cite{etgar-2024, mao-1999, stemple-1968}. 

Finally, unless stated otherwise, all graphs in this paper are connected, and for a positive integer $k$ we use the notation $[k] = \{1,\ldots, k\}$. 
\section{Monotonicity of mutual-visibility sets}
\label{sec:monotonicity}

In this section, for a given visibility set we consider monotonicity of its subsets and monotonicity of its restriction to convex and isometric subgraphs. We recall the previous results and round off the picture so that all four variants are treated systematically. Applying one of our findings we also correct an earlier result from the literature. 

Our starting point is the following result.

\begin{proposition} {\rm \cite[Proposition 2.5]{CiDiDrHeKlYe-2023}}
\label{prop:subset-closed}
If $X$ is a mutual-visibility set (resp.\ outer, total mutual-visibility set) of a graph $G$ and $Y\subseteq X$, then $Y$ is a mutual-visibility set (resp.\ outer, total mutual-visibility set) of $G$.
\end{proposition}

Proposition~\ref{prop:subset-closed} does not hold for dual mutual-visibility sets. For instance, if $x$ and $y$ are adjacent vertices of $C_6$, then $\{x,y\}$ is a $\mud$-set of $C_6$, but neither $\{x\}$ nor $\{y\}$ is a dual mutual-visibility set. 

Dual mutual-visibility therefore stands out because in contrast to the other three types of visibility sets, they are not necessarily closed for taking subsets. On the other hand, all four concepts are monotone for subsets in the following sense. 

\begin{proposition} {\rm \cite[Lemma 5.4]{cicerone-2024b}}
\label{prop:vertex-deleted-subgraphs}
If $X$ is a  mutual-visibility  set (resp.\ outer, dual, or total mutual visibility set) of a graph $G$ and $x\in X$, then $X\setminus \{x\}$ is a  mutual-visibility set (resp.\ outer, dual, or total mutual visibility set) of $G - x$.
\end{proposition}

In the seminal paper on the mutual-visibility, the following useful property was observed. 

\begin{lemma}\label{lem:mu-convex}{\rm \cite[Lemma 2.1]{distefano-2022}}
Let $H$ be a convex subgraph of $G$ and let $X$ be a mutual-visibility set of $G$. Then $X\cap V(H)$ is a mutual-visibility set of $H$.
\end{lemma}

We now show that Lemma~\ref{lem:mu-convex} extends to all the other three mutual-visibility concepts. 

\begin{lemma}\label{lem:three-conveies}
Let $X$ be a dual (outer, total) mutual-visibility set of $G$. If $H$ is a convex subgraph of $G$, then $X\cap V(H)$ is a dual (outer, total) mutual-visibility set of $H$.
\end{lemma}

\proof
Let $X\subseteq V(G)$ and $Y = X\cap V(H)$.

Assume first that $X$ is a dual mutual-visibility set of $G$. We claim that $Y$ is a dual mutual-visibility set of $H$. By Lemma~\ref{lem:mu-convex}, $Y$ is a mutual-visibility set of $H$, hence any two vertices from $Y$ are $Y$-visible in $H$. Consider two vertices $u$ and $v$ from $V(H)\setminus Y$. In $G$, there exists a shortest $u,v$-path $P$ with all internal vertices from $V(G)\setminus X$. Since $H$ is a convex subgraph of $G$, the path $P$ lies completely in $H$. As $V(H)\setminus Y = V(H)\setminus X$, the vertices $u$ and $v$ are $Y$-visible in $H$. We can conclude that $Y$ is a dual-mutual-visibility set of $H$. 

If $X$ is an outer mutual-visibility set of $G$, then, using Lemma~\ref{lem:mu-convex} again, we can proceed as above to prove that $Y$ is an outer mutual-visibility set of $H$. Finally, if $X$ is a total mutual-visibility set of $G$, then combining the above arguments we get that $Y$ is a total mutual-visibility set of $H$.
\qed

Let $G_n$, $n\ge 2$, be the graph obtained from $n$ disjoint $5$-cycles by selecting one edge in each of them and identifying these $n$ edges into a single edge $uv$. Note that $\deg_{G_n}(u) = \deg_{G_n}(v) = n+1$ while the other vertices have degree $2$. In~\cite[Proposition~5.1]{CiDiDrHeKlYe-2023} it was stated that that $\mud(G_n) = n+1$. We now apply Lemma~\ref{lem:three-conveies} to show that this is not the case. The correct result reads as follows.

\begin{proposition}\label{prop:correction}
If $n\ge 2$, then $\mud(G_n) = 2$. 
\end{proposition} 

\proof
Let $X$ be a dual mutual visibility set of $G_n$. Note first that a dual mutual-visibility set of $C_5$ is either the empty set or consists of two adjacent vertices. Since each $5$-cycle of $G_n$ is convex, Lemma~\ref{lem:three-conveies} implies that $X$ restricted to an arbitrary $5$-cycle of $G_n$ is either empty or contains two adjacent vertices.  

Let the vertices of the $i^{\rm th}$ cycle of $G_n$, $i\in [n]$, be $u, x_i, y_i, z_i, v$. Assuming that $X\ne \emptyset$, by the above argument, at least one of the $5$-cycles of $G_n$ has exactly two vertices in $X$. We may assume without loss of generality that this is the cycle $C: u, x_1, y_1, z_1, v$. Up to symmetry, there are three cases to be considered. 

Assume first that $X\cap V(C) = \{u,v\}$. Since $u$ and $v$ lie in each of the $5$-cycles, the above argument yields that $X$ cannot contain further vertices. We may observe that this case is not possible since then $x_1$ and $x_2$ cannot see each other. Assume next that $X\cap V(C) = \{u,x_1\}$. Then the cycle $u, x_2, y_2, z_2, v$ must contain another vertex of $X$ which is adjacent to $u$, and this can only be $x_2$. But then $x_1$ and $x_2$ both belong to $X$ and are not $X$-visible, hence this case is also not possible.  The last case to be considered is $X\cap V(C) = \{x_1, y_1\}$. There is nothing to show if $X$ has no vertices in the other $5$-cycles, hence assume that, without loss of generality, $X\cap \{u, x_2, y_2, z_2, v\} \ne \emptyset$. As $u,v\notin X$, we either have $x_2, y_2\in X$ or $y_2, z_2\in X$. In the first case $x_2$ and $y_1$ are not $X$-visible, in the second case $y_1$ and $y_2$ are not $X$-visible. We can conclude that if $X$ is nonempty, then $X$ intersects only one $5$-cycle. 

To complete the argument we claim that $X = \{x_1, y_1\}$ is a dual mutual-visibility set. Clearly, $x_1$ and $y_1$ are $X$-visible. Consider next arbitrary vertices $x,y\in V(G_n)\setminus \{x_1,y_1\}$. If $x$ and $y$ lie on the same $5$-cycle, they are $X$-visible. And if $x$ and $y$ lie on different $5$-cycles, then  every shortest path between them lies completely inside $V(G_n)\setminus X$. 
\qed

Lemma~\ref{lem:three-conveies} is no longer true if instead of the convexity of the subgraph $H$ we assume that $H$ is isometric. Consider $K_{n,n}$, $n\ge 4$. Then it is not difficult to see that $\mu(K_{n,n}) = \muo(K_{n,n}) = \mud(K_{n,n}) = \mut(K_{n,n}) = 2(n-1)$, and that every largest mutual-visibility set $X$ is of the form $X = V(K_{n,n})\setminus \{u,v\}$, where $u$ and $v$ belong to different bipartition sets of $K_{n,n}$. The subgraph $H = K_{n,n}\setminus \{u,v\} \cong K_{n-1,n-1}$ is isometric, but $X\cap V(H) = V(H)$ is clearly not a mutual-visibility set of $H$ (and hence neither an outer, a dual, or a total mutual-visibility set).  

We also emphasize that the ``converse'' of Lemma~\ref{lem:three-conveies} does not hold. That is, if some set of vertices has the required visibility property on a convex subgraph, it is not always extendable to a set having the same property in the whole graph. For instance, in $C_7$, two adjacent vertices form a convex subgraph and its vertices are of course a total/dual/outer mutual-visibility set of this subgraph. However, two adjacent vertices of $C_7$ do not lie together in a  total, a dual, or an outer \MV\ set. 

We now turn to isometric subgraphs. Note that two adjacent vertices of $C_n$, $n\ge 7$, form a \MV\ set, but the remaining subgraph is not isometric. Similarly,  two antipodal vertices $x$ and $x'$ of $C_n$, $n\ge 6$, form a $\muo$-set of $C_n$, but the graph $C_n - \{x,x'\}$ is not even connected. On the other hand, we have the following positive result. 

\begin{proposition}\label{prop:G-X}
Let $G$ be a connected graph. If $X\subseteq V(G)$ is a dual or a total mutual-visibility set of $G$, then the subgraph $G-X$ is isometric.  
\end{proposition} 

\proof
Assume that $X$ is a dual mutual-visibility set of $G$ and consider any two vertices $x$ and $y$ from $V(G)\setminus X$. Since $X$ is a dual mutual-visibility set, the vertices $x$ and $y$ are $X$-visible, say via a $x,y$-path $P$. But then the path $P$ is also a shortest $x,y$-path in $G-X$, which already implies that $G-X$ is isometric. The same argument applies if $X$ is a total mutual-visibility set. 
\qed

Proposition~\ref{prop:G-X} cannot be strengthened by replacing ``isometric'' with ``convex.'' For instance, if $x$ is a vertex of $C_4$, then $\{x\}$ is a total mutual-visibility set (and hence also a dual mutual-visibility set), but $C_4-x$ is not convex. 

\section{Visibility polynomials}
\label{sec:polynomials}

If $G$ is a graph and $X\subseteq V(G)$, then $X$ is a {\em general position set}~\cite{chandran-2016, manuel-2018} if for any two vertices of $X$, no shortest path between them contains an internal vertex from $X$. In order to better understand these sets, the general position polynomial was introduced in~\cite{irsic-2024}. Here we extend this idea to mutual-visibility sets and pose: 

\begin{definition}
The \emph{visibility polynomial} of a graph $G$ is the polynomial 
$$ \name (G) = \sum_{i \ge 0}r_ix^i\,,$$ 
where $r_i$ is the number of distinct mutual-visibility sets of $G$ with cardinality $i$. 
\end{definition}

Clearly, the degree of $ \name (G)$ is $\mu(G)$, and its constant term is $1$. For instance, if $n\ge 1$, then 
$$\name(P_n) = 1 + nx + \binom{n}{2}x^2\,.$$
In a completely analogous way we define the {\em dual visibility polynomial}, the {\em outer visibility polynomial}, and the {\em total visibility polynomial}, which are, for a given graph $G$, respectively denoted by $\named(G)$, $\nameo(G)$,  {\rm and}\ $\namet(G)$. For paths $P_n$, $n\ge 3$, we have
\begin{align*}
\named(P_n) & = 1 + 2x + 3x^2\,, \\    
\nameo(P_n) & = 1 + nx + x^2\,, \\    
\namet(P_n) & = 1 + 2x + x^2\,. \\    
\end{align*}

As a further example, we determine these four polynomials for balanced complete bipartite graphs. Note that the polynomials for a general complete bipartite graph $K_{m,n}$ can be obtained in the same way but by considering more cases. Here we restrict our attention to the simpler case of $K_{n,n}$.
\begin{proposition} \label{prop:poly-examples}
    For $n \ge 3$, the complete bipartite graph $K_{n,n}$ has the following polynomials:
    \begin{align*}
\name(K_{n,n})\,\, & =
((x+1)^n -x^n)^2+ 2nx^{n+1} + 2x^n\,, \\  
\nameo(K_{n,n}) & = ((x+1)^n -x^n)^2+ 2x^n\,, \\    
\named(K_{n,n}) & = \namet(K_{n,n}) =  ((x+1)^n -x^n)^2\,. \\ 
\end{align*}
    \end{proposition}
    
\proof
    Let $A$ and $B$ be the partite classes of $K_{n,n}$ and consider a set $X \subseteq V(K_{n,n})$. It can be readily checked that $X$ is a mutual-visibility set in each of the following cases:
    \begin{itemize}
        \item[$(a)$] $|X \cap A| \le n-1$ and $|X \cap B| \le n-1$;
        \item[$(b)$] $X=A$ or $X=B$;
        \item[$(c)$] $A \subseteq X$ and $|X \cap B|=1$;
        \item[$(d)$] $B \subseteq X$ and $|X \cap A|=1$.
    \end{itemize}
  Further, if neither of $(a)$-$(d)$ holds, then $X$ contains all vertices from one partite class and at least two vertices, say $u$ and $v$, from the other class. Then $u$ and $v$ are not $X$-visible. Thus $X$ is a mutual-visibility set in $K_{n,n}$ if and only if $X$ satisfies one of $(a)$-$(d)$. This in particular implies that $\mu(K_{n,n}) = 2n-2$.

  If $0 \leq i \leq n+1$, then each $i$-element subset of the vertex set satisfies one of $(a)$-$(d)$ and therefore, we have $r_i= {2n \choose i}$ for the visibility polynomial. If $n+2 \leq i \leq 2n-2$, then only case (a) can be satisfied. There are ${n \choose 2n-i}$ sets $A'$ of cardinality $i$ such that $A\subseteq A'$, and there are ${n \choose 2n-i}$ sets $B'$ of cardinality $i$ such that $B\subseteq B'$. Consequently, $r_i= {2n \choose i} - 2{n \choose 2n-i}$, which in turn implies that 
 \begin{align*}
     \name(K_{n,n})\,\, & =
 \sum_{i=0}^{2n-2} {2n \choose i} x^i - 2\sum_{i=n+2}^{2n-2} {n \choose 2n-i} x^i\\
 & = \sum_{i=0}^{2n} {2n \choose i} x^i - 2x^n\sum_{j=0}^{n} {n \choose j} x^j -2n x^{2n-1} - x^{2n} +2x^n+2nx^{n+1} +2nx^{2n-1} +2x^{2n}\\
 & = (x+1)^{2n} - 2x^n(x+1)^n + x^{2n} +2x^n+2nx^{n+1}\\
 & = ((x+1)^n -x^n)^2 +2nx^{n+1} +2x^n.
 \end{align*}

 For the remaining part of the statement, we note that $X \subseteq V(K_{n,n})$ is an outer mutual-visibility set if and only if condition $(a)$ or $(b)$ holds; and $X$ is a dual mutual-visibility set (or a total mutual-visibility set) if and only if $(a)$ holds. Then, respectively subtracting $2nx^{n+1}$ and $2nx^{n+1} +2x^n$ from $\name(K_{n,n})$ we obtain the polynomials $\nameo(K_{n,n})$ and  $\named(K_{n,n})=\namet(K_{n,n})$.
\qed

Below we give two general properties of the polynomials $\name$, $\nameo$, and $\namet$. For a real number $x$ and an integer $k$ with $x \ge k >0$, the binomial coefficient ${x \choose k}$ is defined as 
\[
{x \choose k} =\displaystyle \prod_{s=1}^k \frac{x-s+1}{s}.
\]
The ``shadow theorem'' of Kruskal~\cite{kruskal} and Katona~\cite{katona} was reformulated by Lov\' asz in~\cite{lovasz} as follows:
\begin{theorem}{\rm \cite{katona, kruskal, lovasz}} \label{thm:kklovasz}
Let $\cF$ be a family of $k$-element sets with $|\cF|= {x \choose k}$ for some real number $x \ge k$. Then the number of different $(k-1)$-element sets covered by $\cF$ is at least ${x \choose k-1}$.  \end{theorem}

Theorem~\ref{thm:kklovasz} and Proposition~\ref{prop:subset-closed} imply the following general property of the coefficients in the polynomials $\name(G)$, $\nameo(G)$, and $\namet(G)$.
\begin{proposition} \label{prop:coefficients}
    Let $G$ be a graph and let $\cP \in \{\name,\nameo,\namet\}$. Suppose that $r_i$ and $r_{i-1}$ are the coefficients of $x^i$ and $x^{i-1}$, respectively, in $\cP(G)$. If $r_i= {z \choose i}$ for a real number $z$, then $r_{i-1} \ge {z \choose i-1}$. 
\end{proposition}

The second general property of the polynomials $\name(G)$, $\nameo(G)$, and $\namet(G)$ is that they can be deduced from the set of all maximal visibility sets as follows, where we set ${\cal P}(X) = (1+x)^{|X|}$ for $X\subseteq V(G)$ and  ${\cal P}\in \{\name, \nameo, \namet\}$.
    
\begin{proposition}
\label{prop:inclusion-exclusion}
Let $G$ be a graph and let ${\cal P}\in \{\name, \nameo, \namet\}$. If $X_1,\ldots, X_n$ is the set of maximal mutual-visibility (resp.\ outer mutual-visibility, resp.\ total mutual-visibility) sets of $G$, then  
$${\cal P}(G) = \sum_{k = 1}^n  (-1)^{k-1} \sum_{\{i_1,\ldots,i_k\}\subseteq [n]} {\cal P}(X_{i_1} \cap \cdots \cap X_{i_k})\,.$$
\end{proposition}

\proof
By Proposition~\ref{prop:subset-closed}, any subset of a mutual-visibility set $X$ is a mutual-visibility set. Hence the contribution of $X$ to $\name(G)$ is $(1+x)^{|X|}$. The formula for $\name(G)$ then follows by the inclusion-exclusion principle. The same argument applies to $\nameo(G)$ and to $\namet(G)$. 
\qed

As an example of the use of Proposition~\ref{prop:inclusion-exclusion}, we will determine  $\nameo(P)$, where $P$ is the Petersen graph. We first infer the following. 

\begin{proposition}
\label{prop:petersen}
Let $P$ be the Petersen graph and $X\subseteq V(P)$. Then $X$ is an outer mutual-visibility set of $P$ if and only if $X$ is an independent set of $P$.     
\end{proposition}

\proof
Assume that $X$ is an outer mutual-visibility set. If two vertices $x$ and $y$ from $X$ are adjacent, and $z$ is a neighbor of $y$ different from $x$, then $x$ and $z$ are not $X$-visible as $P$ is geodetic, a contradiction. 

Conversely, assume that $X$ is an independent set of $P$. If $x,y\in X$, then they are clearly $X$-visible. Assume now that $x\in X$ and $y\notin X$. There is nothing to show if $xy\in E(P)$. Assume hence that $d_P(x,y)=2$ and let $z$ be the common neighbor of $x$ and $y$. Since $X$ is independent and $x\in X$, we have $z\notin X$. We can conclude that a vertex $x\in X$ and a vertex $y\notin X$ are also $X$-visible. 
\qed

Concerning Proposition~\ref{prop:petersen} we remark that one direction of it is a consequence of~\cite[Lemma~5.2]{cicerone-2024b}  which asserts that in a graph of girth at least $5$ every outer mutual-visibility set is an independent set.     

Consider the usual drawing of $P$ and let $u_0, u_1, u_2, u_3, u_4$ be the consecutive vertices of its outer $5$-cycle, and $v_0, v_1, v_2, v_3, v_4$ their respective neighbors in the inner $5$-cycle. Using Proposition~\ref{prop:petersen} and the fact that the independence number of $P$ is $4$, it is straightforward to establish that the only $\muo$-sets of $P$ are: 
$$
\{u_0,u_2,v_3,v_4\}, \{u_1,u_3,v_4,v_0\}, 
\{u_2,u_4,v_0,v_1\}, \{u_3,u_0,v_1,v_2\}, 
\{u_4,u_1,v_2,v_3\}\,.
$$
In addition, there are precisely ten maximal outer mutual-visibility sets of $P$ of size $3$, they are:
\begin{align*}
& \{u_0,v_2,v_3\}, \{u_1,v_3,v_4\},\{u_2,v_4,v_0\},\{u_3,v_0,v_1\},\{u_4,v_1,v_2\}, \\
& \{v_0,u_1,u_4\}, \{v_1,u_2,u_0\},\{v_2,u_3,u_1\},\{v_3,u_4,u_2\},\{v_4,u_0,u_3\}\,.    
\end{align*}
From here, by applying Proposition~\ref{prop:inclusion-exclusion}, we get: 
$$\nameo(P) = 1 + 10x + 30x^2 + 30x^3 + 5x^4\,,$$
where the coefficient at $x^3$ was obtained with computer support. 

We next determine the other three polynomials of $P$. For $\name(P)$, we first state the following result which is of independent interest. 

\begin{proposition}
\label{prop:gp=mu-in-geodetic}
Let $G$ be a geodetic graph and $X\subseteq V(G)$. Then $X$ is a general position set if and only if $X$ is a mutual-visibility set.
\end{proposition}

\proof
A general position set is a mutual-visibility set in general. Hence assume that $X$ is a mutual-visibility set and let $x,y\in X$. Then there exists a shortest $x,y$-path $R$ such that all internal vertices of $R$ lie in $V(G)\setminus X$. But since $G$ is geodetic, $R$ is the unique shortest $u,v$-path, hence $x$ and $y$ lie in general position. We can conclude that $X$ is a general position set. 
\qed    

Since $P$ is geodetic, Proposition~\ref{prop:gp=mu-in-geodetic}  implies that 
$$\name(P) = \psi(P) = 1 + 10x + 45x^2 + 90x^3 + 80x^4 + 30x^5+5x^6\,,$$ 
where $\psi(P)$ is the general position polynomial of $P$. The latter polynomial was introduced in~\cite{irsic-2024}, where the second above equality was also deduced. 

Finally, since $\mud(P) = \mut(P) = 0$, we have $\named(P) = \namet(P) = 1$. 

\section{Gaps in the dual visibility spectrum}
\label{sec:dual-spectrum}

As observed in Section~\ref{sec:monotonicity}, a subset of a dual mutual-visibility set is not necessarily a dual mutual-visibility set. Further, there are graphs admitting $k$-element dual mutual-visibility sets but no $(k-1)$-element ones. For this phenomenon, $C_5$ is the smallest example. We have $\mud(C_5) = 2$, but no single vertex forms a dual mutual-visibility set. This leads to the following concept.
\smallskip

The \emph{dual visibility spectrum} of a graph $G$ is the vector $(r_0, \dots , r_k)$, where $k= \mud(G)$, and $r_i$ is the number of different dual mutual-visibility sets of size $i$ in $G$. Equivalently, the entries $r_0, \dots, r_k$  are the coefficients of $x^0, \dots, x^k$, respectively, in $\named(G)$. We have already observed that $r_0=1$ for every graph.

For example, we have the following dual visibility spectra for cycles:
\begin{itemize}
\item $(1,3,3,1)$ for $C_3$; 
\item $(1,4,4,4)$ for $C_4$; 
\item $(1,0,n)$ for $C_n$ if $n \in \{5,6\}$; and 
\item $(1)$ for $C_n$ with $n \geq 7$.
\end{itemize}

In this section, we show that there can be arbitrarily large gaps, that is, arbitrary zero sequences between positive entries in the dual visibility spectrum of a graph. Moreover, the next result quite surprisingly shows that the spectrum entries can be arbitrarily prescribed, that is, if $r_0=1$, $r_k >0$, and the other entries are arbitrary nonnegative integers, a graph with the given dual visibility spectrum exists.    

\begin{theorem} \label{thm:dual-spectrum}
    For every $k \ge 0$ and nonnegative integers $r_0=1, r_1, \dots, r_k$ with $r_k > 0$, there exists a graph $G$ such that $\mud(G)=k$ and the dual visibility spectrum of $G$ is $(1, r_1, \dots, r_k)$.
\end{theorem}
\proof First, we construct graphs with dual visibility spectra $(1, 0, \dots, 0, 1)$ and $(1, \ell)$, then we build a graph with the spectrum $(1, r_1, \dots , r_k)$.

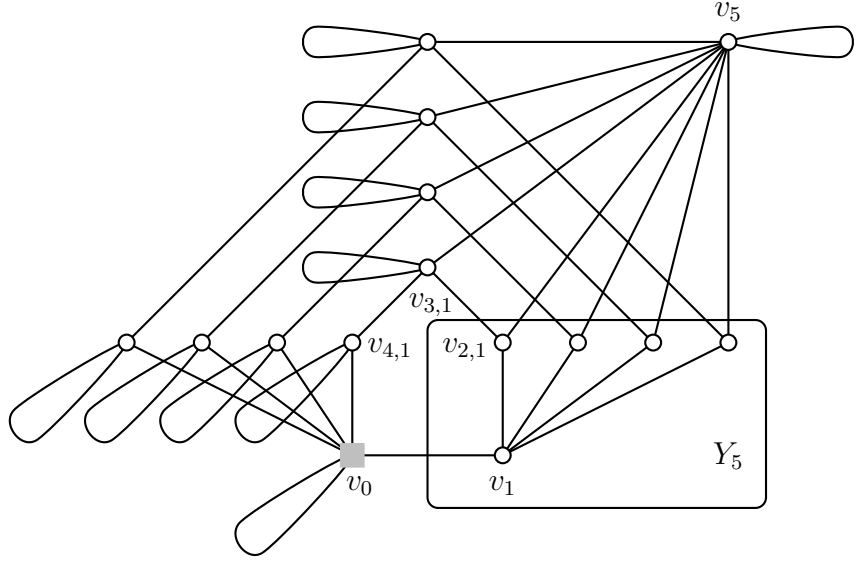
\begin{figure}[ht!]
\begin{center}
\begin{tikzpicture}[scale=1.0,style=thick,x=1cm,y=1cm]
\def\vr{3pt}
\coordinate(v0) at (0.0,0.0);
\coordinate(v1) at (2.0,0.0);
\coordinate(v21) at (2.0,1.5);
\coordinate(v22) at (3.0,1.5);
\coordinate(v23) at (4.0,1.5);
\coordinate(v24) at (5.0,1.5);
\coordinate(v31) at (1.0,2.5);
\coordinate(v32) at (1.0,3.5);
\coordinate(v33) at (1.0,4.5);
\coordinate(v34) at (1.0,5.5);
\coordinate(v41) at (0.0,1.5);
\coordinate(v42) at (-1.0,1.5);
\coordinate(v43) at (-2.0,1.5);
\coordinate(v44) at (-3.0,1.5);
\coordinate(v5) at (5.0,5.5);
\draw (v0) -- (v1);  
\draw (v1) -- (v21) -- (v31) -- (v41) -- (v0);  
\draw (v1) -- (v22) -- (v32) -- (v42) -- (v0);  
\draw (v1) -- (v23) -- (v33) -- (v43) -- (v0);  
\draw (v1) -- (v24) -- (v34) -- (v44) -- (v0);  
\draw (v31) -- (v5) -- (v21);
\draw (v32) -- (v5) -- (v22);
\draw (v33) -- (v5) -- (v23);
\draw (v34) -- (v5) -- (v24);
\draw plot [smooth cycle] coordinates {(0,0) (-1.2,-1.3) (-1.5,-0.9)};
\draw plot [smooth cycle] coordinates {(0,1.5) (-1.2,0.2) (-1.5,0.6)};
\draw plot [smooth cycle] coordinates {(-1,1.5) (-2.2,0.2) (-2.5,0.6)};
\draw plot [smooth cycle] coordinates {(-2,1.5) (-3.2,0.2) (-3.5,0.6)};
\draw plot [smooth cycle] coordinates {(-3,1.5) (-4.2,0.2) (-4.5,0.6)};
\draw plot [smooth cycle] coordinates {(1,2.5) (-0.5,2.7) (-0.5,2.3)};
\draw plot [smooth cycle] coordinates {(1,3.5) (-0.5,3.7) (-0.5,3.3)};
\draw plot [smooth cycle] coordinates {(1,4.5) (-0.5,4.7) (-0.5,4.3)};
\draw plot [smooth cycle] coordinates {(1,5.5) (-0.5,5.7) (-0.5,5.3)};
\draw plot [smooth cycle] coordinates {(5,5.5) (6.5,5.7) (6.5,5.3)};
\draw[rounded corners](1,-0.7)--(1,1.8)--(5.5,1.8)--(5.5,-0.7)--cycle;
\draw(v1)[fill=white] circle(\vr);
\draw(v5)[fill=white] circle(\vr);
\foreach \i in {1,2,3,4} 
{
\draw(v2\i)[fill=white] circle(\vr);
\draw(v3\i)[fill=white] circle(\vr);
\draw(v4\i)[fill=white] circle(\vr);
}
\filldraw[draw=black,color=lightgray] (-0.15,-0.15) rectangle (0.15,0.15);
\node at (0.1,-0.4) {$v_0$};
\node at (2.0,-0.4) {$v_1$};
\node at (5,0) {$Y_5$};
\node at (5,5.9) {$v_5$};
\node at (1.5,1.4) {$v_{2,1}$};
\node at (1.05,2.0) {$v_{3,1}$};
\node at (0.5,1.4) {$v_{4,1}$};
\end{tikzpicture}
\caption{The graph $F_5$}
\label{fig:F_5}
\end{center}
\end{figure}

\paragraph{Construction of $F_t$.}
For every $t \ge 2$, we take $t-1$ $5$-cycles that share the edge $v_0v_1$. For every $i \in [t-1]$, let this cycle be $v_0v_1v_{2,i}v_{3,i}v_{4,i}v_0$. Further, we add a vertex $v_5$ and edges $v_5v_{2,i}$, $v_5v_{3,i}$ for every $i \in  [t-1]$. Let $Y_t= \{v_{2,i}: i \in [t-1]\} \cup \{v_1\}$. To finish the construction, we put a $7$-cycle onto every vertex outside $Y_t$; that is, for each of the vertices $v_0, v_5, v_{3,i}, v_{4,i}$,  where $i \in [t-1]$, we take six new vertices and form a $7$-cycle together with the vertex itself. Vertex $v_0$ is designated as the connecting vertex in $F_t$. The construction is illustrated in Fig.~\ref{fig:F_5} for the case $t=5$, where the gray square emphasizes that $v_0$ is the connecting vertex and where the $7$-cycles are shown as closed ovals. 

Remark that vertex $v_5$ and the incident $7$-cycle may be removed from the graph if $t=2$. In general, some of the $7$-cycles can also be omitted from the construction such that Claim~\ref{claim:1} remains true.

\begin{claim} \label{claim:1}
  For every $t \ge 2$, it holds that $\mud(F_t)=t$ and the only nonempty dual mutual-visibility set of $F_t$ is $Y_t$.
\end{claim}

\proof Suppose that $X$ is a dual mutual-visibility set in $F_t$. Observe that the $7$-cycles in $F_t$ are all convex subgraphs and that $\mud(C_7)=0$. It follows by Lemma~\ref{lem:three-conveies} that $X$ contains no vertices from these $7$-cycles. In particular, $X \subseteq Y_t$. Observe that each of the $t-1$ $5$-cycles is also a convex subgraph in $F_t$. A nonempty dual mutual-visibility set of a $5$-cycle consists of two adjacent vertices. Since $X \subseteq Y_t$, this pair of adjacent vertices may only be $v_1$ and $v_{2,i}$. Therefore, if $v_1 \in X$, then $v_{2,i} \in X$ for all $i \in [t-1]$; and if a vertex $v_{2,i}$ belongs to $X$, we get the same conclusion. We may conclude that $X= \emptyset$ or $X= Y_t$. It can be checked directly that $Y_t$ is a dual mutual-visibility set in $F_t$.  \smallqed

\paragraph{Construction of $F_{1, \ell}$.}
For every $\ell \ge 1$, take an $\ell$-star with a center $v_0$ and leaves  $v_1, \dots , v_\ell$. Then, for every two indices $1\le i < j \le \ell $, add a vertex $u_{i,j}$ and the edges $u_{i,j} v_i$ and $u_{i,j} v_j$. We add edges $u_{i,j} u_{i', j'}$ for every pair of vertices with $1\le i < j \le \ell $ and $1\le i' < j' \le \ell $, that is, the vertices $u_{i,j}$, $1\le i < j \le \ell$, induce a complete subgraph of $F_{1, \ell}$. We put a $7$-cycle onto each vertex outside $Y_{1, \ell} = \{v_1, \ldots, v_\ell\}$ to finish the construction. Vertex $v_0$ is designated as the connecting vertex in $F_{1, \ell}$. We note that $F_{1,1}$ is obtained from $P_2$ by putting a $7$-cycle onto one vertex of it; $F_{1,2}$ can be described as a $4$-cycle where $7$-cycles are put onto two opposite vertices of the $4$-cycle. The construction is illustrated in Fig.~\ref{fig:F_{1,4}} for the case $\ell = 4$, where again the gray square emphasizes that $v_0$ is the connecting vertex and the $7$-cycles are shown as closed ovals. 

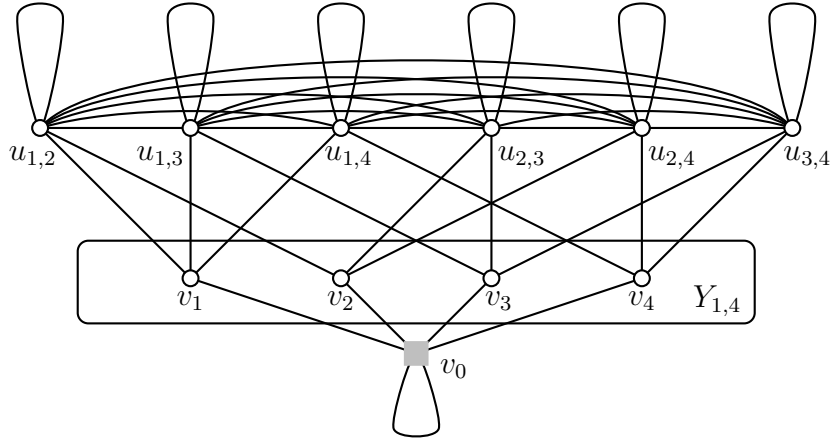
\begin{figure}[ht!]
\begin{center}
\begin{tikzpicture}[scale=1.0,style=thick,x=1cm,y=1cm]
\def\vr{3pt}
\coordinate(v0) at (0.0,0.0);
\coordinate(v1) at (-3.0,1.0);
\coordinate(v2) at (-1.0,1.0);
\coordinate(v3) at (1.0,1.0);
\coordinate(v4) at (3.0,1.0);
\coordinate(v12) at (-5.0,3.0);
\coordinate(v13) at (-3.0,3.0);
\coordinate(v14) at (-1.0,3.0);
\coordinate(v23) at (1.0,3.0);
\coordinate(v24) at (3.0,3.0);
\coordinate(v34) at (5.0,3.0);
\draw (v0) -- (v1);  
\draw (v0) -- (v2);  
\draw (v0) -- (v3);  
\draw (v0) -- (v4);  
\draw (v1) -- (v12) -- (v2);  
\draw (v1) -- (v13) -- (v3);  
\draw (v1) -- (v14) -- (v4);  
\draw (v2) -- (v23) -- (v3);  
\draw (v2) -- (v24) -- (v4);  
\draw (v3) -- (v34) -- (v4);  
\draw (v12) -- (v34);  
\draw (v12) .. controls (-4,3.3) and (-2,3.3) .. (v14);
\draw (v13) .. controls (-2,3.3) and (0,3.3) .. (v23);
\draw (v14) .. controls (0,3.3) and (2,3.3) .. (v24);
\draw (v23) .. controls (2,3.3) and (4,3.3) .. (v34);
\draw (v12) .. controls (-4,3.6) and (0,3.6) .. (v23);
\draw (v13) .. controls (-2,3.6) and (2,3.6) .. (v24);
\draw (v14) .. controls (0,3.6) and (4,3.6) .. (v34);
\draw (v12) .. controls (-4,3.9) and (2,3.9) .. (v24);
\draw (v13) .. controls (-2,3.9) and (4,3.9) .. (v34);
\draw (v12) .. controls (-4,4.2) and (4,4.2) .. (v34);
\draw plot [smooth cycle] coordinates {(-5.0,3.0) (-5.3,4.5) (-4.7,4.5)};
\draw plot [smooth cycle] coordinates {(-3.0,3.0) (-3.3,4.5) (-2.7,4.5)};
\draw plot [smooth cycle] coordinates {(-1.0,3.0) (-1.3,4.5) (-0.7,4.5)};
\draw plot [smooth cycle] coordinates {(1.0,3.0) (0.7,4.5) (1.3,4.5)};
\draw plot [smooth cycle] coordinates {(3.0,3.0) (2.7,4.5) (3.3,4.5)};
\draw plot [smooth cycle] coordinates {(5.0,3.0) (4.7,4.5) (5.3,4.5)};
\draw plot [smooth cycle] coordinates {(0.0,0.0) (-0.3,-1.0) (0.3,-1.0)};
\draw[rounded corners](-4.5,0.4)--(4.5,0.4)--(4.5,1.5)--(-4.5,1.5)--cycle;
\foreach \i in {1,2,3,4} 
{
\draw(v\i)[fill=white] circle(\vr);
}
\draw(v12)[fill=white] circle(\vr);
\draw(v13)[fill=white] circle(\vr);
\draw(v14)[fill=white] circle(\vr);
\draw(v23)[fill=white] circle(\vr);
\draw(v24)[fill=white] circle(\vr);
\draw(v34)[fill=white] circle(\vr);
\filldraw[draw=black,color=lightgray] (-0.15,-0.15) rectangle (0.15,0.15);
\node at (0.5,-0.2) {$v_0$};
\node at (4.0,0.7) {$Y_{1,4}$};
\node at (-3,0.7) {$v_1$};
\node at (-1,0.7) {$v_2$};
\node at (1.1,0.7) {$v_3$};
\node at (3,0.7) {$v_4$};
\node at (-5.1,2.6) {$u_{1,2}$};
\node at (-3.4,2.6) {$u_{1,3}$};
\node at (-0.9,2.6) {$u_{1,4}$};
\node at (1.4,2.6) {$u_{2,3}$};
\node at (3.4,2.6) {$u_{2,4}$};
\node at (5.2,2.6) {$u_{3,4}$};
\end{tikzpicture}
\caption{The graph $F_{1,4}$}
\label{fig:F_{1,4}}
\end{center}
\end{figure}

\begin{claim}  \label{claim:2}
  $\mud(F_{1, \ell})=1$ and there are exactly $\ell$ different $\mud$-sets of $F_{1, \ell}$, namely the sets $\{v_1\}, \dots, \{v_\ell\}$.
\end{claim}
\proof Let $X$ be a nonempty dual mutual-visibility set in $F_{1,\ell}$. Since the $7$-cycles are convex subgraphs in $F_{1,\ell}$ and $\mud(C_7)=0$, we may infer $X \subseteq Y_{1,\ell}$. Suppose now that $|X| \ge 2$. Then at least two different vertices $v_i$ and $v_j$ from $Y_{1,\ell}$ belong to $X$. Since $v_i$ and $v_j$ are the only common neighbors of the (nonadjacent) vertices $v_0$ and $u_{i,j}$, the latter two vertices are not $X$-visible that is a contradiction as both $v_0$ and $u_{i,j}$ are outside $X$. It shows that $|X| =1$. Finally, it suffices to check directly that $X= \{v_i\}$ is a dual mutual-visibility set for every $i \in [\ell]$.  \smallqed

\paragraph{Construction of $F(1, r_1, \dots, r_k)$.}
If $k=0$, we set $F(1)\cong C_7$. If $k \ge 1$, we take the following graphs:
\begin{itemize}
    \item if $r_1 \ge 1$, we take a copy of the graph $F_{1, r_1}$; 
    \item for every $i \ge 2$, if $r_i \ge 1$, we take $r_i$ copies of the graph $F_i$. 
\end{itemize}
The set of these graphs is denoted by $\cG$. Thus, 
$|\cG|= \sum_{i=2}^k r_i$ if $r_1=0$, otherwise $|\cG|= 1+ \sum_{i=2}^k r_i$. Finally, we get $F(1, r_1, \dots, r_k)$, by merging the connecting vertices of the graphs in $\cG$ into one vertex $v^*$.
\begin{claim}  \label{claim:3}
  The dual visibility spectrum of $F(1, r_1, \dots, r_k)$ is $(1, r_1, \dots, r_k)$.
\end{claim}
\proof The statement is true for $k=0$ as the dual visibility spectrum of $C_7$ is $(1)$. From now on, we suppose that $k \ge 1$. Let $G=F(1, r_1, \dots , r_k)$. If $k=1$, then $G\cong F_{1, r_1}$ and the statement follows from Claim~\ref{claim:2}. If $k \ge 2$, $r_k=1$, and $r_i =0$ for all  $i \in [k-1]$, then $G\cong F_{k}$ and the statement follows from Claim~\ref{claim:1}. In the remaining cases, $G$ is constructed from at least two graphs. 

Let $X$ be a dual mutual-visibility set of $G$. Observe that each graph $G_s \in \cG$ is a convex subgraph of $G$ (with the connecting vertex of $G_s$ renamed as $v^*$). By Lemma~\ref{lem:mu-convex}, the set $X\cap V(G_s)$ is a dual mutual-visibility set in $G_s$. 

Suppose that $X$ contains vertices from two different graphs $G_s \in \cG$ and $G_p \in \cG$. If $G_s \cong F_t$ and $G_p \cong F_{t'}$ then, by Claim~\ref{claim:1}, the sets $X \cap V(G_s)$ and $X \cap V(G_p)$ correspond to the dual mutual-visibility sets $Y_t$ and $Y_{t'}$ in $F_t$ and $F_{t'}$. Naming the vertices as in the construction, we consider vertex $v_{1}$ from $G_s$ and $v_{2,1}$ from $G_p$. Both vertices belong to $X$, and the unique shortest path between them goes through the vertex $v_1$ from $G_p$. As the latter vertex is also included in $X$, the set $X$ cannot be a dual mutual-visibility set. In the other case, $G_s \cong F_{1,\ell}$ and $G_p \cong F_{t'}$. Here, we choose vertex $v_i$ from $X\cap V(G_s)$ and consider the shortest path between $v_i$ from $G_s$ and $v_{2,1}$ from $G_p$. The contradiction then comes from the fact that the shortest path is unique and contains vertex $v_1$ from $X\cap V(G_s)$.

We conclude that $X$ cannot intersect two different graphs from $\cG$, and therefore, $X= \emptyset$ or $X$ is a nonempty dual mutual-visibility set in a graph $G_s \in \cG$. By construction of $G$ and by Claims~\ref{claim:1} and~\ref{claim:2}, the dual visibility spectrum of $G$ is  $(1, r_1, \dots, r_k)$. \smallqed
\medskip

 Claim~\ref{claim:3} directly implies the theorem.
 \qed
 
 The following result is a direct corollary of Theorem~\ref{thm:dual-spectrum}.

 \begin{corollary}
     Every polynomial with nonnegative integer coefficients and with a constant term $r_0=1$ is a dual visibility polynomial of some graph.
 \end{corollary}

 By definition, every total mutual-visibility set is a dual mutual-visibility set. All subsets of a $\mut$-set are total and, consequently, dual mutual-visibility sets according to Proposition~\ref{prop:subset-closed}. This establishes the following statement:

\begin{observation}
If $(1,r_1, \dots , r_k)$ is the dual visibility spectrum of a graph $G$ and $i\in [\mut(G)]$, then $r_i \ge {\mut(G) \choose i}$. In particular, there are no gaps in the dual visibility spectrum until the entry $r_j$ with $j=\mut(G)$.
\end{observation}

 We point out a further relation between dual and total mutual-visibility sets.
 \begin{proposition} \label{prop:r1-zero}
     Let $(1,r_1, \dots , r_k)$ be the dual visibility spectrum of a graph $G$. Then $r_1 =0$ if and only if $\mut(G)=0$.
 \end{proposition}
 
\proof
If $\mut(G) >0$, there is a one-element total mutual-visibility set and, by definition, it is also a dual mutual-visibility set. Therefore, we have $r_1 >0$.

If $r_1 >0$, let $X=\{x\}$ be a dual mutual-visibility set. To show that $X$ is also a total mutual-visibility set, we observe that, for every $v \in V(G) \setminus X$, a shortest $x,v$-path never contains an internal vertex from $X$.  It implies $\mut(G) \ge 1$. 
\qed

With respect to the last result we add that the graphs $G$ with $\mut(G) = 0$ were characterized in a different way in~\cite{tian-2024}.

\section{Revisiting total mutual-visibility}
\label{sec:total}

In this section, we characterize total mutual-visibility sets, graphs $G$ with $\mu_t(G) = 1$, and sets which are not total mutual-visibility sets, yet every proper subset is such. 

The vertex $v$ of $G$ is {\em simplicial} if $N_G(v)$ induces a complete subgraph of $G$. The set of simplicial vertices of $G$ is denoted by $S(G)$ and its cardinality by $s(G)$.

To start, let us show the following result. 

\begin{proposition}\label{pro:geodetic graph}
If $G$ is a geodetic graph, then $\mut(G) = s(G)$ and $S(G)$ is the unique $\mut$-set of $G$. 
\end{proposition}

\proof
$S(G)$ is a total mutual-visibility set of $G$ because a vertex from $S(G)$ cannot be an inner vertex of a shortest path. 

To prove that $\mut(G)\leq s(G)$, suppose on the contrary that there exists some $\mut$-set $X$ of $G$ with $|X|\geq s(G) + 1$. Then $X$ contains a vertex $x\notin S(G)$. Let $x_1$ and $x_2$ be two neighbors of $x$ such that $x_1x_2\notin E(G)$. As $x\in X$, and $X$ is a total mutual-visibility set, there exists a vertex $x'\ne x$ such that $x'\in N_G(x_1)\cap N_G(x_2)$. But then there exists at least two shortest $x_1,x_2$-paths, a contradiction.  

We have thus seen that $\mut(G) = s(G)$. Moreover, the above argument also implies that $S(G)$ is the unique $\mut$-set of $G$. 
\qed

In the proof of Proposition~\ref{pro:geodetic graph} it was sufficient to consider only vertices at distance $2$. The announced characterization of total mutual-visibility sets says it is no coincidence. 

\begin{theorem}
\label{thm:distance two}
If $G$ is a connected graph and $X\subseteq V(G)$, then the following statements are equivalent. 
\begin{enumerate}
    \item[(i)] $X$ is a total mutual-visibility set of $G$.
    \item[(ii)] Any two vertices $u$ and $v$ of $G$ with $d_G(u,v) = 2$ are $X$-visible.  
    \item[(iii)] Any two vertices $u$ and $v$ of $G$ with $d_G(u,v) = 2$ satisfy $N_G(u) \cap N_G(v) \not\subseteq X$.
\end{enumerate}
\end{theorem} 

\proof
Let $X\subseteq V(G)$ be a total mutual-visibility set. Then in particular each pair of vertices at distance $2$ is $X$-visible, that is, (i) implies (ii). 

To see that (ii) implies (iii), let $u$ and $v$ be vertices with $d_G(u,v) = 2$. Then by (ii), there exists a shortest $u,v$-path such that its middle vertex, say $w$, does not lie in $X$. Hence $w\in (N_G(u)\cap N_G(v))\setminus X$. 

It remains to prove that (iii) implies (i). That is, we need to show that if (iii) holds, then any two vertices $u',v'\in V(G)$ are $X$-visible. To do so, we proceed by induction on $k = d_G(u',v')$. There is nothing to prove if $k=1$, while if $k=2$, the condition (iii) immediately implies that $u'$ and $v'$ are $X$-visible. Assume now that $k\ge 3$. Let $P$ be a shortest $u',v'$-path, and let $u'=x_0,x_1,x_2$ be its first three vertices. Then $d_G(x_0,x_2) = 2$, hence by (iii) there exists a vertex $y\in N_G(x_0)\cap N_G(x_2)$ such that $y\notin X$. (It is possible that $y=x_1$.) Since $d_G(y,v') = k-1$, the vertices $y$ and $v'$ are $X$-visible by induction. Let $Q$ be a shortest $y,v'$-path such that no internal vertex lies in $X$. Since $y\notin X$, the concatenation of the edge $u'y$ with $Q$ is a shortest $u',v'$-path such that no internal vertex lies in $X$. Hence $u'$ and $v'$ are $X$-visible.  
\qed 

The equivalence between (i) and (iii) in Theorem~\ref{thm:distance two} has been earlier established in~\cite[Theorem 2.3]{Bujtas} for the case of Hamming graphs. 

As already mentioned, in~\cite{tian-2024} the graphs $G$ with $\mut(G) = 0$ were characterized. Moreover, an open problem to characterize the graphs with $\mut(G)=1$ was also posed. In the second main result of this section we solve the problem as follows. For its formulation we recall that a vertex $v$ of a graph $G$ is a {\em bypass vertex}~\cite{tian-2024} if $v$ is not the central vertex of a convex path on three vertices. The number of bypass vertices of $G$ is denoted by $\bp(G)$. 

\begin{theorem} 
\label{thm:mut-1}
For a graph $G$, it holds that $\mut(G)=1$ if and only if $\bp(G) \ge 1$ and every two different bypass vertices $v_1$ and $v_2$ satisfy the following condition:
\begin{itemize}
\item[$(\star)$] there exist nonadjacent vertices $u_1$, $u_2$ with $N_G(u_1) \cap N_G(u_2) = \{v_1, v_2\}$.
\end{itemize}
\end{theorem}
  
\proof
First suppose that $\mut(G)=1$ and $\{v\}$ is a $\mut$-set of $G$. Then $v$ is a bypass vertex and $\bp(G) \ge 1$. Consider now two bypass vertices $v_1$ and $v_2$.  Since $X=\{v_1,v_2\}$ is not a total mutual-visibility set in $G$,  Theorem~\ref{thm:distance two} implies the existence of two vertices $x$ and $y$ with $d_G(x,y)=2$ that satisfy $N_G(x) \cap N_G(y) \subseteq X$. On the other hand, we know the following facts:
\begin{itemize}
    \item $N_G(x) \cap N_G(y) \neq \emptyset$ as $d_G(x,y)=2$;
    \item $N_G(x) \cap N_G(y) \neq \{v_i\}$, for $i \in [2]$, as $v_i $ is a bypass vertex.
    \end{itemize}
Therefore, the only possibility to have $N_G(x) \cap N_G(y) \subseteq X$ is  $N_G(x) \cap N_G(y) = \{v_1, v_2\}$. This proves that every pair of bypass vertices satisfies $(\star)$.
\medskip

To prove the other direction, we take the contrapositive of the implication and assume that $\mut(G)\neq 1$. If $\mut(G)=0$, then $\bp(G) =0$. If $\mut(G) \ge 2$, consider a $2$-element total mutual-visibility set $X=\{ v_1, v_2\}$. Note that, by Proposition~\ref{prop:subset-closed}, such a set exists even if $\mut(G)>2$. Clearly, $v_1$ and $v_2$ must be bypass vertices. We state that $(\star)$ does not hold for $v_1$ and $v_2$. Indeed, the existence of vertices $u_1$ and $u_2$ with $N_G(u_1) \cap N_G(u_2) = \{v_1, v_2\}$ would imply that $u_1$ and $u_2$ are not $X$-visible, a contradiction.
\qed

We note that graphs with $\bp(G) =\ell$ and $\mut(G)=1$ exist for arbitrarily large $\ell$. Graphs $F_{1,\ell}$ constructed in the proof of Theorem~\ref{thm:dual-spectrum} provide such examples for not only $\mud(F_{1,\ell})=1$ but also  $\mut(F_{1,\ell})=1$ holds. 
\medskip

Finally, in view of our considerations in Section~\ref{sec:monotonicity}, and as an application of Theorem~\ref{thm:distance two}, we provide a characterization for sets that are not total-mutual-visibility sets, although all their proper subsets have this property.

\begin{proposition}
    Let $X$ be a nonempty set of vertices in a graph $G$ and suppose that every proper subset $X' \subset X$ is a total \MV\ set in $G$. Then $X$ itself is not a total \MV\ set if and only if there exist two nonadjacent vertices $v_1$ and $v_2$ with $N_G(v_1) \cap N_G(v_2) =X$.
\end{proposition}

\proof
    First observe that $v_1v_2 \notin E(G)$ and $N_G(v_1) \cap N_G(v_2) =X$ imply that each shortest $v_1,v_2$-path contains an internal vertex from $X$. Consequently, $X$ is not a total \MV\ set in $G$.

    Now, assume that $X$ is not a total \MV\ set in $G$, but all proper subsets of $X$ have that property. By Theorem~\ref{thm:distance two}, there exist vertices $v_1$ and $v_2$ with $d_G(v_1,v_2)=2$ that satisfy $N_G(v_1) \cap N_G(v_2) \subseteq X$. Then Theorem~\ref{thm:distance two} also implies that $X'=N_G(v_1) \cap N_G(v_2)$ is not a total \MV\ set of $G$. 
    Therefore, by our condition in the statement, $X'$ is not a proper subset of $X$, and we
   may conclude that $X=X'$, that is, $N_G(v_1) \cap N_G(v_2) =X$ as stated.  
\qed

\section*{Acknowledgments}

We would like to thank Vesna Ir\v{s}i\v{c} Chenoweth for her computer support in calculating the outer visibility polynomial of the Petersen graph. 

This work was supported by the Slovenian Research Agency (ARIS) under the grants P1-0297, N1-0355, and N1-0285.


\begin{thebibliography}{99}

\bibitem{adhikary-2022}
R.~Adhikary, K.~Bose, M.K.~Kundu, B.~Sau,
Mutual visibility on grid by asynchronous luminous robots,
Theoret.\ Comput.\ Sci.\ 922 (2022) 218--247.

\bibitem{axenovich-2024a+}
M.~Axenovich, D.~Liu,
Visibility in hypercubes,
arXiv:2402.04791 [math.CO] (2024).

\bibitem{bollobas-1967}
B.~Bollob\'as,
On a conjecture of {E}rd\H os, {H}ajnal and {M}oon,
Amer.\ Math.\ Monthly 74 (1967) 178--179. 

\bibitem{boruzanli-2024}
G.~Boruzanl\i{} Ekinci, Cs.~Bujt\'as,
Mutual-visibility problems in Kneser and Johnson graphs,
Ars.\ Math.\ Contemp.\ (2024) \url{doi.org/10.26493/1855-3974.3344.4c8}.

\bibitem{BresarYero-2024}
B.~Bre\v{s}ar, I.G.~Yero,
Lower (total) mutual visibility in graphs,
Appl.\ Math.\ Comput.\ 465 (2024) Paper 128411.

\bibitem{Bujtas}
Cs.~Bujt\'{a}s, S.~Klav\v{z}ar, J.~Tian, 
Total mutual-visibility in Hamming graphs,
Opuscula Math. 45 (2025) 63--78.

\bibitem{chandran-2016} 
U.~Chandran S.V., G.J.~Parthasarathy, 
The geodesic irredundant sets in graphs, 
Int.\ J.\ Math.\ Combin.\ 4 (2016) 135--143.

\bibitem{cicerone-2025}
S.~Cicerone, A.~Di Fonso, G.~{Di Stefano}, A.~Navarra, F.~Piselli, Mutual visibility in hypercube-like graphs, 
Lecture Notes Comput.\ Sci.\ 14662 (2024) 192--207.

\bibitem{cicerone-2023a}
S.~Cicerone, G.~{Di Stefano},
Mutual-visibility in distance-hereditary graphs: a linear-time algorithm,
Procedia Comput.\ Sci.\ 223 (2023) 104--111.

\bibitem{CiDiDrHeKlYe-2023}
S.~Cicerone, G.~{Di Stefano}, L.~Dro\v{z}\dj ek, J.~Hed\v{z}et, S.~Klav\v{z}ar, I.G.~Yero,
Variety of mutual-visibility problems in graphs,
Theoret.\ Comput.\ Sci.\ 974 (2023) Paper 114096.

\bibitem{cicerone-2023}
S.~Cicerone, G.~{Di Stefano}, S.~Klav\v{z}ar,
On the mutual-visibility in Cartesian products and in triangle-free graphs,
Appl.\ Math.\ Comput.\ 438 (2023) Paper 127619.

\bibitem{cicerone-2023+}
S.~Cicerone, G.~{Di Stefano}, S.~Klav\v{z}ar, I.G. Yero,
Mutual-visibility in strong products of graphs via total mutual-visibility,
Discrete Appl.\ Math.\ 358 (2024) 136--146. 

\bibitem{cicerone-2024b}
S.~Cicerone, G.~{Di Stefano}, S.~Klav\v{z}ar, I.G.~Yero,
Mutual-visibility problems on graphs of diameter two, 
European J.\ Combin.\ 120 (2024) Paper 103995.

\bibitem{distefano-2022}
G.~{Di Stefano},
Mutual visibility in graphs,
Appl.\ Math. Comput.\ 419 (2022) Paper 126850.

\bibitem{ellis-2024+}
D.~Ellis, M.-R.~Ivan, I.~Leader,
Tur\'an densities for daisies and hypercubes, 
arXiv:2401.16289v6 [math.CO] (2024).

\bibitem{etgar-2024}
A.~Etgar, N.~Linial, 
On the connectivity and diameter of geodetic graphs,
European J.\ Combin.\ 116 (2024) Paper 103886.

\bibitem{irsic-2024}
V.~Ir\v{s}i\v{c}, S.~Klav\v{z}ar, G.~Rus, J.~Tuite,
General position polynomials,
Results Math. 79 (2024) Paper 110. 

\bibitem{katona}
G.O.H.~Katona,
A theorem of finite sets. 
In: Theory of Graphs, Proceedings of Colloquium held at Tihany, Hungary, pp.\ 187--207, Akad\' emia Kiad\' o, Budapest, 1968. 

\bibitem{korze-2024}
D.~Kor\v{z}e, A.~Vesel, 
Mutual-visibility sets in Cartesian products of paths and cycles,
Results Math.\ 79 (2024) Paper 116.

\bibitem{korze-2024+}
D.~Kor\v{z}e, A.~Vesel, 
Variety of mutual-visibility problems in hypercubes,
arXiv:2405.05650  [math.CO] (2024).

\bibitem{kruskal}
J.~Kruskal,
The optimal number of simplices in a complex.
In: Mathematical Optimization Techniques pp.\ 251--278,
Univ.\ California Press, Berkeley, CA,  1963.

\bibitem{kuziak-2023}
D.~Kuziak, J.A.~Rodr\'{\i}guez-Vel\'{a}zquez,
Total mutual-visibility in graphs with emphasis on lexicographic and Cartesian products,
Bull.\ Malays.\ Math.\ Sci.\ Soc.\ 46 (2023) Paper 197.

\bibitem{lovasz}
L. Lov\' asz,
Combinatorial Problems and Exercises. Akad\' emiai Kiad\' o --
North Holland, 1979.
  
\bibitem{manuel-2018}
P.~Manuel, S.~Klav{\v z}ar,
A general position problem in graph theory,
Bull.\ Aust.\ Math.\ Soc.\ 98 (2018) 177--187.

\bibitem{mao-1999}
J.~Mao, D.~Li, 
A new upper bound on the number of edges in a geodetic graph,
Discrete Math.\ 202 (1999) 183--189. 

\bibitem{roy-2025}
D.~Roy, S.~Klav\v{z}ar, A.S.~Lakshmanan,
Mutual-visibility and general position in double graphs and in Mycielskians,
Appl.\ Math.\ Comput.\ 488 (2025) Paper 129131.

\bibitem{stemple-1968}
J.G.~Stemple, M.E.~Watkins, 
On planar geodetic graphs,
J.\ Combin.\ Theory 4 (1968) 101--117.

\bibitem{tian-2024}
J.~Tian, S.~Klav\v{z}ar,
Graphs with total mutual-visibility number zero and total mutual-visibility in Cartesian products,
Discuss.\ Math.\ Graph Theory 44 (2024)  1277--1291.

\bibitem{wessel-1966}
W.~Wessel,
\"Uber eine {K}lasse paarer {G}raphen. {I}. {B}eweis einer {V}ermutung von {E}rd\H os, {H}ajnal und {M}oon,
Wiss.\ Z.\ Tech.\ Hochsch.\ Ilmenau 12 (1966) 253--256. 

\end{thebibliography}
\end{document}